\documentclass[reqno]{amsart}
\usepackage{amssymb}
\usepackage{amsmath}

\makeatletter
\@addtoreset{equation}{section}
\makeatother

\renewcommand\thefigure{\thesection.\@arabic\c@figure}
\renewcommand\thetable{\thesection.\@arabic\c@table}
 
\newtheorem{theorem}{Theorem}[section]
\newtheorem{lemma}[theorem]{Lemma}
\newtheorem{proposition}[theorem]{Proposition}

\newtheorem{remark}[theorem]{Remark}
 
\newcommand{\mc}[1]{{\mathcal #1}}   

\newcommand{\mb}[1]{{\mathbf #1}}
\newcommand{\bb}[1]{{\mathbb #1}}

\begin{document}

\author{M. D. Jara, C. Landim} 

\address{\noindent IMPA, Estrada Dona Castorina 110,
CEP 22460 Rio de Janeiro, Brasil 
\newline
e-mail:  \rm \texttt{monets@impa.br}
}

\address{\noindent IMPA, Estrada Dona Castorina 110,
CEP 22460 Rio de Janeiro, Brasil and CNRS UMR 6085,
Universit\'e de Rouen, 76128 Mont Saint Aignan, France.
\newline
e-mail:  \rm \texttt{landim@impa.br}
}

\title[Nonequilibrium Central Limit Theorem] {Nonequilibrium Central
  Limit Theorem for a Tagged Particle in Symmetric Simple Exclusion} 

\begin{abstract}
  We prove a nonequilibirum central limit theorem for the position of
  a tagged particle in the one-dimensional nearest-neighbor symmetric
  simple exclusion process under diffusive scaling starting from a
  Bernoulli product measure associated to a smooth profile $\rho_0:\bb
  R\to [0,1]$.
\end{abstract}

\maketitle

\section{Introduction}
\label{sec0}

The asymptotic behavior of a tagged particle appears as one of the
central problems in the theory of interacting particle systems and
remains mostly unsolved.

The first important result on the position of a tagged particle in the
diffusive scaling is due to Kipnis and Varadhan \cite{kv}. By proving
an invariance principle for additive functionals of reversible Markov
processes, Kipnis and Varadhan deduced an equilibrium central limit
theorem for the position of a tagged particle in symmetric simple
exclusion processes. This result was extended by Varadhan \cite{v} for
mean-zero asymmetric exclusion processes, through an invariance
principle for Markov processes with generator satisfying a sector
condition; and by Sethuraman, Varadhan and Yau \cite{svy} to
asymmetric exclusion processes in dimension $d\ge 3$, relaxing the
sector condition by a graded sector condition. In these three contexts
the authors prove that
$$
\frac{X_{tN^2} - E[X_{tN^2}]}{N}
$$
converges in law, as $N\uparrow\infty$, to a Brownian motion with 
diffusion coefficient given by a variational formula. Here $X_t$
stands for the position of the tagged particle at time $t$.

The nonequilibrium picture is much less clear. Even a law of large
numbers for a tagged particle starting from a Bernoulli product
measure with slowly varying parameter seems still out of reach.
Rezakhanlou \cite{r} proved a propagation of chaos result which states
that the average behavior of tagged particles is described by
diffusion process. A large deviations from this diffusive limit in
dimension $d\ge 3$ was obtained by Quastel, Rezakhanlou and Varadhan
\cite{qrv} .

We prove in this article the first nonequilibrium central limit
theorem for a tagged particle. Consider the one-dimensional nearest
neighbor symmetric situation. In this context, as already observed by
Arratia \cite{a}, the scaling changes dramatically since to displace
the tagged particle from the origin to a site $N>0$, all particles
between the origin and $N$ need to move to the right of $N$.  This
observation relates the asymptotic behavior of the tagged particle to
the hydrodynamic behavior of the system. The correct scaling for the
law of large numbers should therefore be $X_{t N^2}/N$ and we expect
$(X_{t N^2} - E[X_{t N^2}])/\sqrt{N}$ to converge to a Gaussian
variable.

The central limit theorem in equilibrium was obtained by Rost and
Vares \cite{rv} for a slightly different model. They proved that for
each fixed $t>0$, $X_{t N^2} /\sqrt{N}$ converges to a fractional
Brownian motion $W_t$ with variance given by $E[W_t^2] = \alpha
t^{1/2}$.  We extend their result to the nonequilibrium case.

The idea of the proof is to relate the position of the tagged particle
to the well known hydrodynamic behavior of the symmetric exclusion
process. Since particles cannot jump over other particles, the
position of the tagged particle is determined by the current over one
bond and the density profile of particles. Therefore, a nonequilibrium
central limit theorem for the position of the tagged particle follows
from a joint central limit theorem for the current and the density
profile. Since the current over a bond can itself, at least formally,
be written as the difference between the mass at the right of the bond
at time $t$ and the mass at time $0$, a central limit theorem for the
position of the tagged particle should follow from a nonequilibrium
central limit theorem for the density field.  This is the content of
the article.

There are three main ingredients in the proof. In Section \ref{sec2}
we present a nonequilibrium central limit theorem for the current over
a bond and show how it relates to the fluctuations of the density
field.  In section \ref{sec4} we obtain a formula which relates the
position of the tagged particle to the current over one bond and the
density field. Finally, in Section \ref{sec5} we present a sharp
estimate on the difference of the solution of the hydrodynamic
equation and the solution of a discretized version of the hydrodynamic
equation.

\section{Notation and Results}
\label{sec1}

The nearest neighbor one-dimensional symmetric exclusion process is a
Markov process on $\{0,1\}^{\mathbb Z}$ which can be described as
follows. Particles are initially distributed over $\bb Z$ in such a
way that each site is occupied by at most one particle. A particle at
a site $x$ waits for an exponential time and then jumps to $x \pm 1$
provided the site is vacant.  Otherwise the jump is suppressed and the
process starts again.

The state space of this Markov process is denoted by
$\mathcal{X}=\{0,1\}^\bb Z$ and the configurations by the Greek letter
$\eta$, so that $\eta(x)=1$ if site $x$ is occupied for the
configuration $\eta$ and $0$ otherwise. The generator $L_N$ of the
process speeded up by $N^2$ is given by
\begin{equation*}
(L_N f) (\eta) \;=\;  N^2 \sum_{x \in \bb Z}
[f(\sigma^{x,x+1}\eta)-f(\eta)] \;,
\end{equation*}
where $\sigma^{x,x+1}\eta$ is the configuration obtained from
$\eta$ by interchanging the occupation variables $\eta(x)$ and
$\eta(x+1)$:
$$
(\sigma^{x,x+1}\eta)(z) \;=\;
\left\{
\begin{array}{ll}
\eta(x+1) & \text{if $z=x$,} \\
\eta(x) & \text{if $z=x+1$,} \\
\eta(z) & \text{otherwise.}
\end{array}
\right.
$$

For each configuration $\eta$, denote by $\pi(\eta)$ the positive
measure on $\bb R$ obtained by assigning mass $N^{-1}$ to each
particle:
$$
\pi (\eta) \;=\; N^{-1} \sum_{x \in \bb Z} \eta(x)\delta_{x/N}
$$
and let $\pi_t =\pi(\eta_t)$.

Fix a profile $\rho_0 :\bb R \to [0,1]$ with the first four derivatives limited.
Denote by $\nu^N_{\rho_0(\cdot)}$ the product measure on $\mc X$ associated to
$\rho_0$:
$$
\nu^N_{\rho_0(\cdot)} \{\eta, \eta(x)=1\} \;=\; \rho_0(x/N)
$$
for $x$ in $\bb Z$. For each $N\ge 1$ and each measure $\mu$ on
$\mc X$, denote by $\bb P_\mu$ the probability on the path space
$D(\bb R_+, \mc X)$ induced by the measure $\mu$ and the Markov
process with generator $L_N$. Expectation with respect to $\bb P_\mu$
is denoted by $\bb E_\mu$. Note that we omitted the dependence of the
probability $\bb P_\mu$ on $N$ to keep notation simple. This
convention is adopted below for several other quantities which also
depend on $N$. The hydrodynamic behavior of the symmetric simple
exclusion process is trivial and described by the heat equation.

\begin{theorem}
\label{s1}
Fix a profile $\rho_0 :\bb R \to [0,1]$.  Then, for all time $t \geq
0$, under $\bb P_{\nu^N_{\rho_0(\cdot)}}$ the sequence of random
measures $\pi_t$ converges in probability to the absolutely continuous
measure $\rho(t,u) du$ whose density $\rho$ is the solution of the heat
equation with initial condition $\rho_0$:
\begin{equation}
\label{eq:02}
\left\{
\begin{array}{l}
\partial_t\rho = \Delta \rho \\
\rho(0,\cdot ) = \rho_0(\cdot). \\
\end{array}
\right.
\end{equation}
Here and below, $\Delta$ stands for the Laplacian. 
\end{theorem}

This theorem establishes a law of large numbers for the empirical
measure.  To state the central limit theorem some notation is
required. For $k\ge 0$, denote by $\mc H_k$ the Hilbert space induced
by smooth rapidly decreasing functions and the scalar product $<\cdot,
\cdot>_k$ defined by
$$
<f,g>_k \;=\; < f , (x^2-\Delta)^k g>\;,
$$
where $<\cdot, \cdot>$ stands for the usual scalar product in $\bb
R^d$. Notice that $\mc H_0 = L^2(\bb R^d)$ and denote by $\mc H_{-k}$
the dual of $\mc H_k$.

Let $\rho^N_t(x) = \bb E_{\nu^N_{\rho_0(\cdot)}} [\eta_t (x)]$. A
trivial computation shows that $\rho^N_t(x)$ is the solution of the
discrete heat equation:
\begin{equation}
\label{eq:05}
\left\{
\begin{array}{l}
\partial_t \rho^N_t (x) \;=\; \Delta_N \rho^N_t (x)\;, \\
\rho^N_0 (x) = \rho_0(x/N)\;,
\end{array}
\right.
\end{equation}
where $(\Delta_N h)(x) = N^2 \sum_{y,|y-x|=1} [h(y) - h(x)]$.

Fix $k\ge 4$ and denote by $\{Y^N_t, t\ge 0\}$ the so called density
field, a $\mc H_{-k}$-valued process given by
$$
Y_t^N (G) \;=\; \frac{1}{\sqrt{N}} \sum_{x \in \bb Z} 
G(x/N) \{ \eta_t (x) - \rho^N_t(x)\}
$$
for $G$ in $\mc H_k$. Denote by $Q_N$ the probability measure on the
path space $D(\bb R_+, \mc H_{-k})$ induced by the process $Y_t^N$ and
the measure $\nu^N_{\rho_0(\cdot)}$. Next result is due to Galves,
Kipnis and Spohn in dimension 1 and to Ravishankar \cite{rav} in
dimension $d\ge 2$.

\begin{theorem}
\label{s2}
The sequence $Q_N$ converges to $Q$, the probability measure
concentrated on $C(\bb R_+, \mc H_{-k})$ corresponding to the
Orsntein-Uhlenbeck process $Y_t$ with mean zero and covariance given
by
$$
\mathbb{E}[Y_t(H)Y_s(G)] \;=\; \int_{\bb R} (T_{t-s} H) \, G \, \chi_s
\;-\; \int_0^s dr\, \int_{\bb R} (T_{t-r}H) \, (T_{s-r}G) \, 
\{\partial_r \chi_r - \Delta \chi_r\}
$$
for $0 \leq s < t$ and $G,H \in \mathcal H_k$. In this formula, $\{T_t
: t\ge 0\}$ stands for the semigroup associated to the Laplacian and 
$\chi_s$ for the function $\chi(s,u) = \rho(s,u)[1-\rho(s,u)]$.
\end{theorem}

Note that in the case of the heat equation, $\partial_r \chi_r -
\Delta \chi_r = 2(\partial_x \rho)^2$. Also, in the equilibrium case,
$\chi$ is constant in space and time so that  the second term vanishes
and we recover the equilibrium covariances. Finally, integrating by
parts twice the expression with $\Delta \chi_r$, we rewrite the
limiting covariances as
\begin{equation}
\label{eq:12}
\mathbb{E}[Y_t(H)Y_s(G)] \;=\; \int_{\bb R} (T_{t} H) \, (T_s G) \, \chi_0
\;+\; 2 \int_0^s dr\, \int_{\bb R} (\nabla T_{t-r}H) \, (\nabla T_{s-r}G) \, 
\chi_r \;,
\end{equation}
where $\nabla f$ is the space derivative of $f$.

\medskip 

We examine in this article nonequilibrium central limit theorems for
the current through a bond and the position of a tagged particle.  For
a bond $(x,x+1)$, denote by $J_{x,x+1} (t)$ the current over this bond.
This is the total number of jumps from site $x$ to site $x+1$ in the
time interval $[0,t]$ minus the total number of jumps from site $x+1$
to site $x$ in the same time interval. 

\begin{theorem}
\label{s3}
Fix $u$ in $\bb R$ and let 
$$
Z^N_t\;=\; \frac{1}{\sqrt{N}} \Big\{ J_{x_N, x_N+1}(t) - 
E_{\nu^N_{\rho_0(\cdot)}} [J_{x_N, x_N+1}(t)]\Big\}\;,
$$
where $x_N = [uN]$. Then, for every $k\ge 1$ and every $0\le t_1
<\cdots <t_k$, $(Z_{t_1}^N , \dots, Z_{t_k}^N)$ converges in law to a
Gaussian vector $(Z_{t_1}, \dots, Z_{t_k})$ with covariance given by
\begin{eqnarray*}
E[Z_s Z_t] &=&
\int_{-\infty}^0 dx \,  P[B_s\le x] \, P[B_t\le x] \, \chi_0(x) \\
&+& \int_0^{\infty} dx \,  P[B_s\ge x] \, P[B_t\ge x] \, \chi_0(x) \\
&+& 2 \int_0^s dr \int_{-\infty}^{\infty} dx \, p_{t-r}(0,x) 
\, p_{s-r}(0,x) \, \chi_r(x)
\end{eqnarray*}
provided $s\le t$ and $u=0$. In this formula, $B_t$ is a standard
Brownian motion starting from the origin and $p_t(x,y)$ is the
Gaussian kernel.
\end{theorem}

By translation invariance, in the case $u\not = 0$, we just need to
translate $\chi$ by $-u$ in the covariance.

Let $H_0 = \mb 1\{[0, \infty)\}$. The covariance appearing in the
previous theorem is easy to understand. Formally the current $N^{-1/2}
J_{-1,0}(t)$ centered by its mean corresponds to $Y^N_t(H_0)$ $-
Y^N_0(H_0)$ since both processes increase (resp. decrease) by
$N^{-1/2}$ whenever a particle jumps from $-1$ to $0$ (resp. $0$ to
$-1$).  The limiting covariance $E[Z_s Z_t]$ corresponds to the formal
covariance
$$
E\Big[ \big\{ Y_t (H_0) - Y_0(H_0) \big\} 
\big\{ Y_s (H_0) - Y_0(H_0) \big\} \Big]\; .
$$

Denote by $\nu^{N,*}_{\rho_0(\cdot)}$ the measure
$\nu^N_{\rho_0(\cdot)}$ conditioned to have a particle at the origin.

\begin{remark}
\label{s8}
The law of large numbers and the central limit theorem for the
empirical measure and for the current starting from
$\nu^{N,*}_{\rho_0(\cdot)}$ follow from the law of large numbers and
the central limit theorem for the empirical measure and the current
starting from the measure $\nu^N_{\rho_0(\cdot)}$ since we may couple
both processes in such a way that they differ at most at one site at
any given time.
\end{remark}

Fix a profile $\rho_0$ with the first four derivatives limited, and consider
the product measure $\nu^{N,*}_{\rho_0(\cdot)}$. Denote by $X_t$ the position
at time $t\ge 0$ of the particle initially at the origin. A law of large
numbers for $X_t$ follows from the hydrodynamic behavior of the
process:

\begin{theorem}
\label{s4}
Fix $t\ge 0$.  $X_t/N$ converges in $\bb
P_{\nu^{N,*}_{\rho_0(\cdot)}}$-probability to $u_t$, the solution of
$$
\dot u_t \;=\; - \frac{(\partial_u \rho)(t,u_t)}{\rho(t,u_t)}\;\cdot
$$
\end{theorem}

Note that the solution of the previous equation is given by
$$
\int^{u_t}_0 du\, \rho(t,u) \;=\; -\int^t_0 ds\, (\partial_u\rho) 
(s,0) \;.
$$

\begin{theorem} 
\label{s7}
Assume that $\rho_0$ has a bounded fourth derivative. Let $W_t =
N^{-1/2}$ $(X_t - N u_t)$. Under $\bb P_{\nu^{N,*}_{\rho_0(\cdot)}}$, For
every $k\ge 1$ and every $0\le t_1 <\cdots <t_k$, $(W_{t_1}^N , \dots,
W_{t_k}^N)$ converges in law to a Gaussian vector $(W_{t_1}, \dots,
W_{t_k})$ with covariance given by
\begin{eqnarray*}
\rho(s,u_s) \rho(t,u_t) E[W_s W_t] &=& 
\int_{-\infty}^0 dx \,  P_{u_s}[B_s\le x] 
\, P_{u_t}[B_t\le x] \, \chi_0(x) \\
&+& \int_0^{\infty} dx \,  P_{u_s} [B_s\ge x] \, P_{u_t}[B_t\ge x] 
\, \chi_0(x) \\
&+& 2 \int_0^s dr \int_{-\infty}^{\infty} dx \, p_{t-r}(u_t,x) 
\, p_{s-r}(u_s,x) \, \chi_r(x)\;.
\end{eqnarray*}
In this formula, $P_u$ stands for the probability corresponding to a
standard Brownian motion starting from $u$.
\end{theorem}

The assumption made on the smoothness of $\rho_0$ appears because in
the proof of Theorem \ref{s7} we need a sharp estimate on the
difference of the discrete approximation of the heat equation
(\ref{eq:05}) and the heat equation (\ref{eq:02}). In section
\ref{sec5} we show that there exists a finite constant $C_0$ for which
$|\rho^N_t(x) - \rho(t,x/N)|\le C_0 t N^{-2}$ for all $N\ge 1$, $x$ in
$\bb Z$ and $t\ge 0$ under the assumption that $\rho_0$ has a bounded
fourth derivative.

\section{Nonequilibrium fluctuations of the Current}
\label{sec2}

Suppose for a moment that the profile $\rho_0$ has a compact support.
Then, $\eta_0$ is almost surely a configuration with a finite number
of particles, and it is easy to see that we have a simple formula for
the current $J_{-1,0}(t)$:
\begin{equation}
\label{eq1}
J_{-1,0}(t) = \sum_{x \geq 0} \eta_t(x)-\eta_0(x) \;.
\end{equation}
In particular, we can write $J_{-1,0}(t)$ in terms of the fluctuation
field:
$$
\frac{1}{\sqrt{N}} \left\{J_{-1,0} (t) - E_{\nu^N_{\rho_0(\cdot)}}
[J_{-1,0}(t)]\right\} \;=\;  Y^N_t(H_0) - Y^N_0(H_0\;),
$$
where $H_a$ is the indicator function of the interval $[a,\infty)$:
$$
H_a (u) \;=\; \mb 1\{[a,\infty)\} (u)\;.
$$
Since the profile has compact support, it is possible to define
$Y_t(H_0)$ as the limit $Y_t(G_n)$ for some sequence $G_n$ of compact
supported function converging to $H_0$ on compact subsets of $\bb R$
and to prove that $Y^N_t(H_0)$, defined in a similar way, converges to
$Y_t(H_0)$. 

In the general case, however, when $\rho_0$ is an arbitrary profile,
neither formula (\ref{eq1}) makes sense, nor the fluctuation field
$Y_t^N(H_0)$ is well defined.  Nevertheless, there is a way to
calculate the fluctuations of the current by appropriated
approximations of the function $G$, as made by Rost and Vares
\cite{rv} in the equilibrium case.

Define the sequence $\{G_n : n\ge 1\}$ of approximating functions of
$H_0$ by 
$$
G_n(u) = \{1 - (u/n)\}^+ \mb 1\{u\ge 0\}\;.
$$
From here we use the next convention: if $X$ is a random variable, we
denote by $\overline{X}$ the centered variable $X-
E_{\nu^N_{\rho_0(\cdot)}}[X]$. 

\begin{proposition}
\label{s5}
For every $t\ge 0$,
$$
\lim_{n \to \infty}
\bb E_{\nu^N_{\rho_0(\cdot)}} \Big[ N^{-1/2} \overline{J}_{-1,0}(t)
- Y^N_t(G_n) + Y^N_0(G_n) \Big]^2 \;=\; 0
$$
uniformly in N.
\end{proposition}

\begin{proof}
Clearly,
$$
M_{x,x+1}(t) \; :=\;  J_{x,x+1} (t)\;-\; N^2
\int_0^t ds\, \{\eta_s (x) -\eta_s(x+1)\}
$$
is a martingale with quadratic variation given by
$$
<M_{x,x+1}>_t \;=\; N^2 \int_0^t ds \,
\{\eta_s (x) -\eta_s(x+1)\}^2 \;.
$$
The goal is to express the difference $Y^N_t(G_n) - Y^N_0(G_n)$ in
terms of the martingales $M_{x,x+1}(t)$ and to notice that these
martingales are orthogonal, since they have no common jumps.

Since
$$
J_{x-1,x}(t) - J_{x,x+1}(t) \;=\; \eta_t(x)-\eta_0(x)
$$
for all $x$ in $\bb Z^d$, $t\ge 0$,
$$
Y^N_t(G_n) - Y^N_0(G_n) \;=\; N^{-1/2} \sum_{x \in \bb Z}
G_n (x/N) \{ \overline{J}_{x-1,x}(t) - \overline{J}_{x,x+1}(t)\}\;.
$$
A summation by parts and the explicit form of $G_n$ permits to rewrite
this expression as
$$
N^{-1/2} \overline{J}_{-1,0}(t) 
- N^{-1/2} \sum_{x = 1}^{nN} \frac{1}{nN}
\overline{J}_{x-1,x}(t)\;.
$$
Representing the currents $J_{x,x+1} (t)$ in terms of the martingales
$M_{x,x+1}(t)$, we obtain that
\begin{eqnarray*}
\!\!\!\!\!\!\!\!\!\!\!\!\!\!\!\! &&
N^{-1/2} \overline{J}_0(t) - \left[Y^N_t(G_n)-Y^N_0(G_n)\right] \\
\!\!\!\!\!\!\!\!\!\!\!\!\!\!\!\! && \quad
=\;  \frac{1}{\sqrt{N}}\sum_{x=1}^{nN} \frac{1}{nN} M_{x-1,x}(t)
\;+\; \frac{1}{\sqrt{N}} \int_0^t ds\, \frac{N}{n}
\left[ \overline{\eta}_s(0) -
\overline{\eta}_s(nN) \right] \;.
\end{eqnarray*}

We claim that the martingale and the integral term converge to $0$ in
${L}^2(\bb{P}_{\nu^N_{\rho_0(\cdot)}})$. In fact, since the
martingales are orthogonal, estimating their quadratic variations by
$tN^2$, an elementary computation shows that
\begin{equation*}
\bb E_{\nu^N_{\rho_0(\cdot)}} \Big [\frac{1}{\sqrt{N}}\sum_{x=1}^{nN}
\frac{1}{nN} M_{x-1,x}(t) \Big]^2 \; \leq\; \frac{t}{n}\;.
\end{equation*}

The integral term is more demanding, because in non-equilibrium the
two-point correlations are not easy to estimate. Expanding the square
we have that
\begin{eqnarray*}
\!\!\!\!\!\!\!\!\!\!\!\!\!\!\!\! &&
\bb E_{\nu^N_{\rho_0(\cdot)}} \Big[\frac{1}{\sqrt{N}} \int_0^t ds\,
\frac{N}{n} [\overline{\eta}_s(0) -\overline{\eta}_s(nN) ] ds \Big]^2 \\
\!\!\!\!\!\!\!\!\!\!\!\!\!\!\!\! && \quad
=\; \frac{2N}{n^2} \int^t_0 ds \int^s_0 dr \,
\bb E_{\nu^N_{\rho_0(\cdot)}} \Big[ \big(\overline{\eta}_s(0)
-\overline{\eta}_s(nN)\big) \big(\overline{\eta}_r(0)
- \overline{\eta}_r(nN)\big) \Big]\; .
\end{eqnarray*}
By Lemma \ref{s6} the previous expression is less than or equal to
$C_0 t^{5/2} n^{-2}$ for some finite constant $C_0$ depending only on
$\rho_0$. This concludes the proof of the proposition.
\end{proof}

A central limit theorem for the current $\overline{J}_{-1,0}(t)$ is a
consequence of this proposition.

\noindent{\bf Proof of Theorem \ref{s3}}.
Fix $t\ge 0$ and $n\ge 1$. By approximating $G_n$ in $L^2(\bb R) \cap
L^1(\bb R)$ by a sequence $\{H_{n,k} : k\ge 1\}$ of smooth functions
with compact support, recalling Theorem \ref{s2}, we show that
$Y^N_t(G_n)$ converges in law to a Gaussian variable denoted by
$Y_t(G_n)$.

By Proposition \ref{s5}, $\{Y^N_t(G_n) - Y^N_0(G_n) : n\ge 1\}$ is a
Cauchy sequence uniformly in $N$. In particular, $Y_t(G_n) - Y_0(G_n)$
is a Cauchy sequence and converges to a Gaussian limit denoted by $Y_t
(H_0) - Y_0(H_0)$. Therefore, by Proposition \ref{s5},
$\overline{J}_{-1,0}(t)$ converges in law to $Y_t (H_0) - Y_0(H_0)$.

The same argument show that any vector $(\overline{J}_{-1,0}(t_1),
\dots, \overline{J}_{-1,0}(t_k))$ converges in law to $(Y_{t_1}
(H_0) - Y_0(H_0), \dots, Y_{t_k} (H_0) - Y_0(H_0))$. The covariances
can be computed since by (\ref{eq:12})
\begin{eqnarray*}
\!\!\!\!\!\!\!\!\!\!\!\!\!\!\! &&
E\Big[ \big\{ Y_{t} (H_0) - Y_0(H_0) \big\} \big\{ Y_{s} (H_0) -
Y_0(H_0) \big\}  \Big] \\
\!\!\!\!\!\!\!\!\!\!\!\!\!\!\! && \quad
=\;\lim_{n\to\infty} E\Big[ \big\{ Y_{t} (G_n) - Y_0(G_n) \big\} 
\big\{ Y_{s} (G_n) - Y_0(G_n) \big\}  \Big] \\
\!\!\!\!\!\!\!\!\!\!\!\!\!\!\! && \quad
=\; \lim_{n\to\infty} \Big\{ \int_{\bb R} \big\{ (T_{t} G_n) (T_s G_n)  +
G_n^2  - (T_t G_n) G_n  - (T_s G_n) G_n \big\} \chi_0 \\
\!\!\!\!\!\!\!\!\!\!\!\!\!\!\! && \qquad\qquad\qquad
\;+\; 2 \int_0^s dr \int_{\bb R} (\nabla T_{t-r} G_n) \, 
(\nabla  T_{s-r} G_n) \chi_r \Big\}\;.
\end{eqnarray*} 
A long but elementary computation permits to recover the expression
presented in the statement of the theorem. This concludes the
proof. \qed

\medskip

We conclude this section with some elementary estimates on two points
correlation functions. For $0\le s\le t$ and $x\neq y$ in $\bb Z$, let
$$
\phi(t; x,y) = E_{\nu^N_{\rho_0(\cdot)}} [\eta_t(x) ; \eta_t(y)]
\;,\quad
\phi(s,t; x,y) = E_{\nu^N_{\rho_0(\cdot)}} [\eta_s(x) ; \eta_t(y)]\; .
$$
In this formula and below, $E_\mu[f;g]$ stands for the covariance of
$f$ and $g$ with respect to $\mu$.

\begin{lemma}
\label{s6}
There exists a finite constant $C_0=C_0(\rho_0)$ depending only on the
initial profile $\rho_0$ such that
$$
\sup_{x,y \in \bb Z} |\phi(t; x,y)| \;\le\; \frac{C_0 \sqrt{t}}{N}\; ,
\quad
\sup_{x,y \in \bb Z} |\phi(s,t;x,y)| \;\leq\; \frac{C_0}{N} \Big\{
\sqrt{s} + \frac 1{\sqrt{t-s}}\Big\}\;.
$$
\end{lemma}

The first statement is a particular case of an estimate proved in
\cite{fpsv}. In sake of completeness, we present an elementary proof
of this lemma.

\begin{proof}
  Let $L_2$ be the generator of $2$ nearest-neighbor symmetric simple
  exclusion processes on $\bb Z$. An elementary computation shows that
  $\phi(t; x,y)$ satisfies the difference equation
$$
\left\{
\begin{array}{l}
(\partial_t \phi)(t;x,y) = N^2 (L_2 \phi)(t;x,y) - \mb 1\{|x-y|=1\}
N^2 [\rho^N(t,x) - \rho^N(t,y)]^2 \;, \\
\phi(0;x,y) = 0\;.
\end{array}
\right.
$$
This equation has an explicit solution  which is (negative and)
absolutely bounded by
$$
C_0(\rho_0) \int_0^t ds\, \bb P_{x,y} \big[ |X_s - Y_s|=1\big]
$$
for $C_0 = \Vert \partial \rho_0 \Vert_\infty^2$.  In this formula,
$(X_s,Y_s)$ represent the position of the symmetric exclusion process
speeded up by $N^2$ and starting from $\{x,y\}$. A coupling argument
shows that $\bb P_{x,y} [ |X_s - Y_s|=1] \le \bb P_{x,y}^0 [ |X_s -
Y_s|=1]$ where in the second probability particles are evolving
independently.  Since $\bb P_{x,y}^0 [ |X_s - Y_s|=1] \le C
(sN^2)^{-1/2}$, the first part of the lemma is proved.

To prove the second statement, recall that we denote by $\Delta_N$ the
discrete Laplacian in $\bb Z$. $\phi(t;y) = \phi(s,t;x,y)$ satisfies
the difference equation
$$
\left\{
\begin{array}{l}
(\partial_t \phi) (t;y) = (\Delta_N \phi) (t;y) \\
\phi(s;y) = \phi (s;x,y) \text{ if $y\neq x$,} \\
\phi(s;y) = \rho^N(s,x) [1-\rho^N(s,x)] \text{ for $y= x$.}
\end{array}
\right.
$$
This equation has an explicit solution
$$
\phi(s;y) \;=\; \sum_{z\neq x} p_{t-s}(y,z) \phi(s;x,z)
\;+\; p_{t-s}(y,x) \rho^N(s,x) [1-\rho^N(s,x)]\;,
$$
where $p_s(x,y)$ stands for the transition probability of a nearest
neighbor symmetric random walk speeded up by $N^2$. The first part of
the lemma together with well known estimates on $p_s$ permit to conclude.
\end{proof}

\section{Law of Large Numbers for the Tagged Particle}
\label{sec3}

In this section we assume the initial measure to be
$\nu^{N,*}_{\rho_0(\cdot)}$, the product measure
$\nu^N_{\rho_0(\cdot)}$ conditioned to have a particle at the origin.
Keep in mind Remark \ref{s8}.

Fix a positive integer $n$. The tagged particle is at the right of $n$
at time $t$ if and only if the total number of particles in the
interval $\{0, \dots, n-1\}$ is less than or equal to the current
$J_{-1,0}(t)$: 
\begin{equation}
\label{eq:01}
\{ X_t \geq n \} \;=\; 
\{ J_{-1,0} (t) \geq \sum_{x=0}^{n-1} \eta_t(x) \} \;.
\end{equation}
This equation indicates that a law of large numbers and a central
limit theorem for the position of the tagged particle are intimately
connected to the joint asymptotic behavior of the current and the
empirical measure. We prove in this section the law of large numbers.

Denote by $\lceil a \rceil$ the smallest integer larger than or equal
to $a$. Fix $u>0$ and set $n=\lceil uN \rceil$ in (\ref{eq:01}) to
obtain that
\begin{equation}
\label{eq:04}
\{ X_t \geq uN \} \;=\; 
\Big \{ N^{-1} J_{-1,0} (t) \; \geq \;\, <\pi^N_t, \mb 1\{[0,u]\}> 
\;+\; O(N^{-1}) \Big\} \;.
\end{equation}

By Theorem \ref{s1}, $<\pi^N_t, \mb 1\{[0,u]\}>$ converges in probability to 
$\int^u_0 \rho(t,w) dw$, where $\rho$ is the solution of the heat
equation (\ref{eq:02}).

On the other hand, the law of large numbers for $J_{-1,0} (t)$ under
$\bb P_{\nu^N_{\rho_0(\cdot)}}$ is an elementary consequence of the
central limit theorem proved in the last section and the convergence
of the expectation of $N^{-1} J_{-1,0} (t)$. By the martingale
decomposition of the current and by Theorem \ref{s10},
\begin{eqnarray*}
\bb E_{\nu^N_{\rho_0(\cdot)}} \big[ N^{-1} J_{-1,0} (t) \big] &=&
\int_0^t ds \, N \big[ \rho^N_s (-1) - \rho^N_s (0) \big] \\ 
&=& - \int^t_0 \partial_u \rho(s,0) ds \;+\; O(N^{-1})\; .
\end{eqnarray*}
Hence, $N^{-1} J_{-1,0} (t)$ converges in probability to $- \int^t_0
\partial_u \rho(s,0) ds$. 

In view of (\ref{eq:04}) and the law of large numbers for the current
and the empirical measure,
$$
\lim_{N \to \infty} \mathbb{P}_{\nu^{N,*}_{\rho_0(\cdot)}} 
\big [N^{-1} X_t \ge u \big] \;=\;
\begin{cases}
0 & \textrm{ if } -\int^t_0 \partial_u \rho(s,0) ds < \int^u_0
\rho(t,w) dw \;, \\
1 & \textrm{ if } -\int^t_0 \partial_u \rho(s,0) ds > \int^u_0
\rho(t,w) dw\;. \\
\end{cases}
$$

By symmetry around the origin, a similar statement holds for $u<0$.
Thus, $X^N_t/N$ converges to $u_t$ in probability, where $u_t$ is
the solution of the implicit equation
$$
\int_0^{u_t} \rho(t,w) dw \;=\;  -\int_0^t \partial_u\rho_s(0) ds\; .
$$

\section{Central Limit Theorem for the Tagged Particle}
\label{sec4}

In this section we prove Theorem \ref{s7} developing the ideas of the
previous section. Assume first that $u_t> 0$ and fix $a$ in $\bb R$.
By equation (\ref{eq:01}), the set $\{X_t \ge N u_t + a\sqrt{N}\}$ is
equal to the set in which
\begin{equation}
\label{eq:03}
\overline{J}_{-1,0} (t) \geq 
\sum_{x=0}^{Nu_t} \overline{\eta}_t (x)
\;+\; \sum_{x=1}^{a\sqrt{N}  -1} \eta_t(x+Nu_t) 
\;-\;  \Big\{ \mathbb{E}_{\nu^{N}_{\rho_0(\cdot)}} [J_{-1,0}(t)]
-\sum_{x=0}^{Nu_t} \rho^N_t(x)]\Big\} \;,
\end{equation}
where $\rho^N_t(x)$ is the solution of the discrete heat equation
(\ref{eq:05}). 

We claim that second term on the right hand side of (\ref{eq:03})
divided by $\sqrt{N}$ converges to its mean in $L^2$. Indeed, by Lemma
\ref{s6}, its variance is bounded by $C_0 a^2 N^{-1}$ for some finite
constant $C_0$. Since by Theorem \ref{s10},
$$
\frac 1{\sqrt{N}} \sum_{x=1}^{a\sqrt{N}  -1} \rho^N_t(x+Nu_t)
$$
converges to $a \rho(t,u_t)$, the second term on the right hand side of
(\ref{eq:03}) converges in probability to $a \rho(t,u_t)$.

An elementary computation based on the definition of $u_t$ and on
Theorem \ref{s10} shows that the third term on the right hand side of
(\ref{eq:03}) divided by $\sqrt{N}$ vanishes as $N\uparrow\infty$.

Finally, by Proposition \ref{s5}, for fixed $t$, $N^{-1/2} \{
\overline{J}_{-1,0} (t) - \sum_{x=0}^{Nu_t} \overline{\eta}_t (x)\}$
behaves as $Y^N_t(G_n) - Y^N_0(G_n) - Y^N_t(\mb 1\{[0,u_t]\})$, as
$N\uparrow\infty$, $n\uparrow\infty$. Repeating the arguments
presented at the beginning of the proof of Theorem \ref{s3}, we show
that this latter variable converges in law to a centered Gaussian
variable, denoted by $W_t$, and which is formally equal to
$Y_t(H_{u_t})-Y_0(H_0)$.  

Up to this point we proved that
$$
\lim_{N\to\infty} \bb P_{\nu^{N}_{\rho_0(\cdot)}} \Big[
\frac{X_t - u_t N}{\sqrt{N}} \ge  a\Big] \;=\;
P[W_t \ge a \rho(t,u_t)]
$$
provided $u_t>0$. The same arguments permit to prove the same
statement in the case $u_t=0$, $a>0$. By symmetry around the origin,
we can recover the other cases: $u_t<0$ and $a$ in $\bb R$, $u_t=0$
and $a<0$.

Putting all these facts together, we conclude that for each fixed $t$,
$(X_t-N u_t)/\sqrt{N}$ converges in distribution to the Gaussian
$W_t/\rho(t,u_t) = [Y_t(H_{u_t})-Y_0(H_0)]/\rho(t,u_t)$.  The same
arguments show that any vector $(N^{-1/2} [X_{t_1}-N u_{t_1}], \dots,
N^{-1/2} [X_{t_k}-N u_{t_k}])$ converges to the corresponding centered
Gaussian vector.

It remains to compute the covariances, which is long but elementary.

\section{Appendix}
\label{sec5}

In sake of completeness, we present in this section a result on the
approximation of the heat equation by solutions of discrete heat
equations.

Fix a profile $\rho_0 :\bb R\to\bb R$ with a bounded fourth
derivative.  Let $\rho : \bb R_+\times \bb R\to\bb R$ be the solution
of the heat equation with initial profile $\rho_0$:
$$
\left\{
\begin{array}{l}
\partial_t \rho(t,x) = \partial_x^2 \rho(t,x) \\
\rho(0,x) = \rho_0(x)\;.
\end{array}
\right.
$$

Recall that we denote by $\Delta_N$ the discrete Laplacian.  For each
$N \in \mathbb{N}$, define $\rho^N_t(x)$ as the solution of the system
of ordinary differential equations
\begin{equation}
\label{eq:08}
\left\{
\begin{array}{l}
(d/dt) \rho^N_t(x) = (\Delta_N \rho^N_t) (x) \\
\rho^N_0(x) = \rho_0 (x/N)\;.
\end{array}
\right.
\end{equation}

The main result of this section asserts that $\rho^N$ approximates
$\rho$ up to order $N^{-2}$:

\begin{theorem}
\label{s10}
Assume that $\rho_0: \mathbb{R} \to [0,1]$ is a function with a
bounded fourth derivative. There exists a finite constant $C_0$ such
that
$$
\Big| \rho^N_t(x) - \rho(t, \frac{x}{N}) \Big| 
\;\leq\; \frac{C_0t}{N^2}
$$
for all  $N\ge 1$, $t \geq 0$, $x \in \bb Z$.
\end{theorem}

An easy way to prove this statement is to introduce a time discrete
approximation of the heat equation.  For each $N$ in $\mathbb{N}$ and
each $\delta > 0$, we define $\rho_l^{\delta, N}(k)$, $k$ in $\bb Z$,
$l \geq 0$ by the recurrence formula
\begin{equation}
\label{eq:07}
\left\{
\begin{array}{l}
\rho_{l+1}^{\delta, N}(k) = \rho_l^{\delta, N}(k) + 
        \delta N^2 [\rho_l^{\delta, N}(k+1)+\rho_l^{\delta, N}(k-1)
                -2\rho_l^{\delta,N}(k)]\\
\rho_0^{\delta, N}(k) = \rho_0(k/N)\;.
\end{array}
\right.
\end{equation}

We now recall two well known propositions whose combination leads to
the proof of Theorem \ref{s10}.  The first one states that the
solution of (\ref{eq:07}) converges as $\delta\downarrow 0$ to the
solution of (\ref{eq:08}) uniformly on compact sets. The second one
furnishes a bound on the distance between the solution of the discrete
equation (\ref{eq:07}) and the solution of the heat equation.

For $a$ in $\bb R$, denote by $\lfloor a \rfloor$ the largest integer
smaller or equal to $a$.

\begin{proposition}
\label{s9}
For each $N\ge 1$,
$$
\lim_{\delta \to 0} \rho^{N,\delta}_{\lfloor t/\delta \rfloor}(k)
\;=\; \rho^N_t(k)
$$
uniformly on compacts of $\bb R_+\times \bb Z$.
\end{proposition}

\begin{proposition}
\label{s11}
Suppose that $\delta N^2 < 1/2$. Then, there exist a finite constant
$C_0=C_0(\rho_0)$ such that
$$
\Big|\, \rho^{\delta,N}_l(k) - \rho\Big(\delta l , k/N \Big)  \, \Big |
\;\leq\; C_0 \Big\{ \delta^2 l +\frac{\delta l}{N^2} \Big\}
$$
for all $l \geq 0, k \in \bb Z$.
\end{proposition}

Clearly, Theorem \ref{s10} is an immediate consequence of Propositions
\ref{s9} and \ref{s11}. Proposition \ref{s9} is a consequence of
Proposition \ref{s11} and the Cauchy-Peano existance theorem for
ordinary differential equations. Proposition \ref{s11} is a standard
result on numerical analysis (see \cite{t} for instance).

\end{document}